\documentclass{article}

\def\ra{\rightarrow}
\def\ms{\mapsto}

\def\mod{\mathop{\;\rm mod}\nolimits}

\def\Pic{\mathop{\rm Pic}\nolimits}
\def\Card{\mathop{\rm Card}}

\def\Z{\hbox{\bf Z}}

\def\Q{\hbox{\bf Q}}

\def\F{\hbox{\bf F}}

\def\eps{\varepsilon}

\def\ni{\noindent}

\begin{document}

\begin{center}
{\large {\bf Courbes de genre $3$ avec $S_3$ comme groupe
d'automorphismes}}

\vspace{7ex}

{\sc Jean-Fran\c cois Mestre.}
\end{center}

\vspace{10ex}

\section{Introduction}

Soit $k$ un corps, $a_1,a_2,a_3,a_4,X,Y,Z$ des ind\'etermin\'ees, 
$T_1=X+Y+Z$, $T_2=XY+YZ+ZX$, $T_3=XYZ$;
la courbe projective plane $C_{a_1,a_2,a_3,a_4}$ \`a coefficients dans $k(a_1,a_2,a_3,a_4)$ d'\'equation
$$a_1T_1^4+a_2T_1^2T_2+a_3T_1T_3+a_4T_2^2=0$$ 
a  comme groupe d'automorphismes le groupe sym\'etrique $S_3$
agissant de la fa\c con naturelle sur $\{X,Y,Z\}$. 

Nous \'etudions ici cette famille de courbes en vue d'applications au probl\`eme
de la recherche de courbes de genre $g$, d\'efinies sur un corps fini $\F_q$, $q=p^n$, $p$
premier, ayant beaucoup de points. 

Plus pr\'ecis\'ement, notons $N_q(g)$ le maximum du nombre de points de $C(\F_q)$, lorsque $C$ parcourt l'ensemble des courbes de genre $g$ d\'efinies
sur $\F_q$, et soit $m_q$ la partie enti\`ere de $2\sqrt{q}$. La borne de Weil, raffin\'ee par Serre ([Serre Harvard]), donne la majoration $N_q(g)\leq q+1+gm_q$.  Une courbe $C$ d\'efinie
sur $\F_q$ dont le nombre de points sur $\F_q$ est \'egal \`a $q+1+gm_q$ est dite
{\it optimale}, et le {\it d\'efaut d'optimalit\'e} est la quantit\'e $D_q(g)=q+1+gm_q-N_q(g)$.

Notons que, pour $q$ et $g$ donn\'es, il peut ne pas exister de courbe optimale,
y compris dans le cas du genre $1$. \`A ce propos,
Serre a pos\'e la question suivante:

\medskip
\ni
{\it Pour $g$ fix\'e, existe-t-il une constante $c(g)$ telle que, pour tout $q$,
$D_q(g)\leq c(g)$?}

\medskip

Dans le cas du genre $3$, rappelons les r\'esultats suivants, en relation avec ceux
obtenus ici:

\medskip

\begin{itemize}

\item Top ([Top]) a trouv\'e de nombreux cas de courbes de genre $3$ optimales
en faisant une recherche syst\'ematique, en particulier pour $q\leq 100$.

\item Lauter ([Lauter2]) a prouv\'e que, pour tout $q$, il existe une courbe $C$ de
genre $3$ telle que $|q+1-\Card C(\F_q)|\geq 3m_q-3.$

\item Auer et Top ([Auer-Top]) ont prouv\'e que, pour $n$ impair, $D_{3^n}(3)\leq 21$.

\item Lorsque $n$ est pair, et que $p\equiv 3\mod 4$, il existe ([Ibukiyama], p. $2$) une courbe
optimale de genre $3$, qui plus est hyperelliptique. Pour $n$ pair $\geq 4$, et $p=2$,
il existe une courbe de genre $3$ optimale ([Na-Rit2]). Pour $p\equiv 1\mod 4$ et $n\equiv 2\mod 4$,
on sait aussi qu'il existe une courbe optimale.

\item Pour $p=2$, Nart et Ritzenthaler ([Na-Rit1]) ont
prouv\'e qu'il existe une courbe optimale pour une infinit\'e de $n$ impairs.

\end{itemize}

\medskip
Plusieurs de ces r\'esultats  ont \'et\'e obtenus en consid\'erant la famille de courbes
$ax^4+by^4+cz^4+ey^2z^2+fz^2x^2+gx^2y^2=0$, dites {\it quartiques de Ciani} ([Ciani], [La-Ritz]), qui poss\`edent
un groupe d'automorphismes isomorphe \`a $(\Z/2\Z)^2$. 

\medskip
 
L'utilisation des courbes $C_{a_1,a_2,a_3,a_4}$ ci-dessus,  contenant le groupe
sym\'etrique $S_3$ comme groupe d'automorphismes, permet d'obtenir les 
r\'esultats suivants:

\medskip
\ni
{\sc Th\'eor\`eme.}
{\it
\begin{itemize}
\item Soit $q=3^n$, $n$ impair. Si $E$ est une courbe elliptique d\'efinie
sur $\F_q$ dont l'invariant modulaire n'est pas dans $\F_3$, il existe
une courbe de genre $3$ dont la jacobienne est isog\`ene sur $\F_q$ \`a
$E^3$.
\item Soit $q=3^n$,
avec $n$ impair $\geq 7$; si $m_q$ n'est pas divisible
par $3$ (donc pour une infinit\'e de $n$),  il existe une courbe
de genre $3$ optimale (et aussi une courbe dont le nombre de points
est $q+1-3m_q$).

Dans le cas o\`u $n$ est pair,
il existe des courbes optimales non hyperelliptiques. 

\item Soit $q=7^n$, $n$ impair, et $a$ un entier premier \`a $7$, $|a|\leq m_q$; si $a$ est
divisible par $3$, et  est un carr\'e  $\mod 7$,
il existe une courbe
de genre $3$ dont le nombre de points est $q+1-3a$. Si $a\equiv 1,4,5,7,8$ ou $11\mod 12$,
il existe une courbe de genre $3$ dont le nombre de points est $q+1-3a$ et une courbe
de genre $3$ dont le nombre de points est $q+1+3a$.
\end{itemize}
}

\medskip
\ni
{\sc Corollaire.}
{\it
\begin{itemize}

\item Pour tout $n$, on a  $D_{3^n}(3)\leq 3$.

\item Pour tout $n$, on a $D_{7^n}(3)\leq 9$
\end{itemize}
}

Par ailleurs, les courbes $C_{a_1,a_2,a_3,a_4}$ 
permettent fr\'equemment d'obtenir des courbes optimales:
par exemple (cf. $2.2$), pour $p$ premier
$\leq 10000$, 
on obtient, par sp\'ecialisation convenable de courbes de cette famille, une courbe optimale sur $\F_p$ dans environ
$90\%$ des cas.

\medskip

La premi\`ere section est consacr\'ee au cas o\`u la caract\'eristique
de $k$ est diff\'erente de $2$ et $3$;  on y donne d'abord un formulaire concernant les 
courbes $C_{a_1,a_2,a_3,a_4}$. Leur jacobienne est isog\`ene
\`a $E_1^2\times E_2$, o\`u $E_1$ (resp. $E_2$) est la courbe quotient de $C$ par l'un quelconque
des sous-groupes d'ordre $2$ (resp. par le sous-groupe d'ordre $3$) de $S_3$,
et
il existe un isomorphisme galoisien commutant \`a l'accouplement de Weil entre
les points d'exposant $3$ de $E_1$ et de $E_2$. 

R\'eciproquement, si $E_1$ et $E_2$ sont deux courbes elliptiques sur $k$,
d'invariants distincts et distincts de $0$ et $1728$, 
il existe
douze courbes de genre $3$ sur $\overline{k}$ de la forme $C_{a_1,a_2,a_3,a_4}$ dont la jacobienne
est isog\`ene \`a $E_1^2\times E_2$. 

\medskip

Dans la sous-section suivante, on applique ces formules au cas o\`u $E_1$ et $E_2$ 
ont le m\^eme invariant modulaire, ce qui donne une famille de courbes de genre
$3$ \`a un param\`etre $t$. Cela permet en particulier de construire une courbe
de genre $3$ dont la jacobienne est isog\`ene sur $\Q$ au cube  d'une courbe \`a multiplications 
complexes, dans le cas o\`u l'anneau d'endomorphismes de cette courbe a comme discriminant
$-3,-7,-19,-43,-67,-163,-16$ ou $-28$.  

Nous traitons ensuite le cas o\`u $p=7$, avant
de finir dans la derni\`ere  section par le cas de la caract\'eristique $3$.

\section{Caract\'eristique diff\'erente de $2$ et $3$}

\ni
{\sc Notations.} 
\begin{itemize}
\item Soit $j\ne 0,1728$. Dans ce qui suit, on note $E(j)$
la courbe elliptique d'\'equation 
$$y^2=x^3-3\frac{j}{j-1728}x
-\frac{3j}{4(j-1728)}.$$ Cette courbe est d'invariant $j$.

\item  On note $\rho\in \overline{k}$ une racine  de $x^2+x+1$.
Remarquons que $(\rho:\rho^2:1)$ est un point de la courbe
$C_{a_1,a_2,a_3,a_4}$.

\item Si $E$ est une courbe elliptique, on note $E[3]$ le groupe de ses
points d'exposant $3$.

\end{itemize}

\subsection{Isog\'enies d'ordre $3$ des courbes elliptiques}
Soit $E$ une courbe elliptique
d'\'equation
$y^2=p(x)$, avec $p(x)=x^3+ax+b$; l'\'equation aux abscisses de ses points d'ordre $3$
est alors $f_3(x)=2 p p"-p'^2=3x^4+6ax^2+12xb-a^2.$ Par une homographie convenable,
on peut ramener les quatre racines de $f_3$ \`a $\infty,1,\rho,\rho^2$.

L'ensemble des birapports des $24$ permutations de ces
racines est \'egal \`a $\{-\rho,-\rho^2\}$, donc le groupe des homographies conservant les quatre
racines de $f_3$ est isomorphe au groupe altern\'e $A_4$. Si $f:\;E[3]\ra
E[3]$ est un isomorphisme pr\'eservant l'accouplement de Weil, l'\'equation
aux abscisses de $f$ restreinte \`a $E[3]-\{0\}$ est donn\'ee par une
telle homographie.

Par suite, si $E'$ est une seconde courbe elliptique, d'\'equation
$y^2=x^3+a'x+b'$, et si $g_3=3x^4+6a'x^2+12xb'-a'^2$, il existe
douze homographies envoyant les racines de $f_3$ sur celles de $g_3$,
et elles correspondent aux isomorphismes de $E[3]$ sur $E'[3]$ commutant
\`a l'accouplement de Weil.

\medskip
\ni
La courbe elliptique universelle poss\'edant un point d'ordre $3$ est la courbe
d'\'equation $$E_t:\;y^2+xy+ty=x^3,$$
le point $(0,0)$ \'etant d'ordre $3$ et $t$ 'etant un
param\`etre.  Son invariant modulaire est \'egal \`a $$j_3(t)=\frac{(24t-1)^3}{t^3(27t-1)}.$$

La courbe quotient de $E_t$ par le groupe d'ordre $3$ engendr\'e par $(0,0)$ est la 
tordue quadratique par $-3$ de la courbe $E_{-t-1/27}.$

\subsection{Formulaire}
Si $a_3$ ou $a_4$ est nul, la courbe $C_{a_1,a_2,a_3,a_4}$  n'est pas
irr\'eductible; on peut donc supposer $a_4=1$ et $a_3\neq 0$;  de plus,
comme la caract\'eristique de $k$ est ici suppos\'ee 
diff\'erente de $2$ et $3$,
on peut se ramener \`a $a_3=2$ par une transformation de la forme
$X'=X+aT_1,Y'=Y+aT_1,Z'=Z+aT_1$, qui commute aux matrices de permutation et qui est inversible pour $3a+1\neq 0$, en prenant
$a=\displaystyle{\frac{2-a_3}{3a_3}}$.

\ni
On d\'esigne d\'esormais par $C_{a_1,a_2}$ la courbe
d'\'equation
$$a_1 T_1^4+a_2T_1^2T_2+2T_1T_3+T_2^2=0.$$

Son discriminant vaut
$256\, \left( 27\,{a_1}+5+9\,{a_2} \right) d^3$, o\`u

\medskip
\noindent
$ d= -432\,{a_1}
+72\,{{ a_2}}^{2}+76\,{{a_2}}^{3}-216\,{{ a_2}}^{2}{ a_1}+432
\,{{ a_1}}^{2}-432\,{ a_1}\,{t a_2}+27\,{{ a_2}}^{4}
.$

\medskip

On note \'egalement $q_{a_1,a_3}(x,y)$ la forme non homog\`ene de
l'\'equation de $C_{a_1,a_2}$ obtenue en
posant $x=X/Z,y=Y/Z$.

\medskip

\ni
{\sc Proposition $2.1.$ }

{\it 
$1$)  
La courbe $E_1$ quotient de $C_{a_1,a_2}$ par le 
groupe engendr\'e par 
l'involution  $(X,Y,Z)\ms (Y,X,Z)$ est une courbe elliptique 
d'\'equation affine
$$y^2=x^3-(3+2a_2)x -4a_1+2+2a_2+a_2^2$$
et d'invariant modulaire $$J_1=6912\frac{(2a_2+3)^3}{d}.$$

Le morphisme
$\phi_1:\;C_{a_1,a_2}\ra E_1$ est donn\'e par 
$$\phi_1:\;\left\{\begin{array}{ll}
x=&\frac{X+Y-Z}{T_1},\\
y=&\frac{2(Z^2-XY)+a_2T_1^2+4T_2}{T_1^2}\end{array}\right.$$

Le point 
\`a l'infini de $E_1$ est l'image par $\phi_1$ de $(\rho:\rho^2:1)$
et l'on a
$$\phi_1^*(\frac{dx}{y})=
\frac{(y-x)dx}{q_y},$$
o\`u $q_y=\displaystyle{\frac{\partial }{\partial y}q_{a_1,a_2}}$. 

\medskip
Si l'on quotiente $C_{a_1,a_2}$ par le groupe engendr\'e par
l'involution $\tau:\;(X,Y,Z)\ms (X,Z,Y)$ (resp. $(X,Y,Z)\ms (X,Z,Y)$), 
on trouve la m\^eme courbe
$E_1$; l'image r\'eciproque de $\frac{dx}{y}$ est 
$\frac{(x-1)dx}{q_y}$ (resp. 
$\frac{(1-y)dx}{q_y}$); 

La somme de ces trois formes diff\'erentielles vaut $0$, et elles 
sont deux-\`a-deux lin\'eairement ind\'ependantes. 

La jacobienne de $C_{a_1,a_2}$ est donc isog\`ene \`a $E_1^2\times E_2$,
o\`u $E_2$ est une courbe elliptique.

\medskip

$2$) La courbe $E_2$ est la jacobienne de  
de la courbe de genre $1$ quotient de
$C_{a_1,a_2}$ par le groupe d'ordre $3$ engendr\'e par
$\sigma:\;(X,Y,Z)\ms (Y,Z,X)$; une \'equation de $E_2$ est
$$y^2=x^3+Ax+B,$$
avec

$$\begin{array}{ll}

A=&-48-128\,a_{{2}}-648\,{a_{{2}}}^{2}+3456\,a_{{1}}a_{{2}}-648\,{a_{{2}}
}^{3}+1944\,a_{{1}}{a_{{2}}}^{2}\\&-243\,{a_{{2}}}^{4}+3168\,a_{{1}}-3888
\,{a_{{1}}}^{2}\\
&\\

B=&512\,a_{{2}}-2400\,{a_{{2}}}^{2}+46080\,a_{{1}}a_{{2}}-7200\,{a_{{2}}}
^{3}+77760\,a_{{1}}{a_{{2}}}^{2}\\&-9720\,{a_{{2}}}^{4}-124416\,{a_{{1}}}
^{2}a_{{2}}$$ $$+54432\,a_{{1}}{a_{{2}}}^{3}-5832\,{a_{{2}}}^{5}-69984\,{a_
{{1}}}^{2}{a_{{2}}}^{2}\\&+17496\,a_{{1}}{a_{{2}}}^{4}-1458\,{a_{{2}}}^{6
}+128+19328\,a_{{1}}$$ $$-120960\,{a_{{1}}}^{2}+93312\,{a_{{1}}}^{3}
\end{array}$$

$3$) Le quotient $J_1/J_2$ des invariants modulaires de $E_1$ et $E_2$ est un
cube $M_3^3$, o\`u
$$M_3=\frac{A}{16\left( 27\,a_{{1}}+5+9\,a_{{2}}
 \right)  \left( 3+2\,a_{{2}} \right)}.$$

$4$) Il existe un isomorphisme galoisien, commutant \`a 
l'accouplement de Weil, entre $E_1[3]$  
et $E_2[3]$. 

Plus pr\'ecis\'ement, il existe une
homographie d\'efinie sur $k$ envoyant les abscisses des points
d'ordre $3$ de $E_1$ sur les abscisses des points d'ordre $3$ de
$E_2$. 

$5$) Soit $J_1$ un \'el\'ement de $k-\{0,1728\}$. L'ensemble des courbes $C_{a_1,a_2}$
telles que $E_1$ a comme invariant $J_1$ est une famille \`a un param\`etre
$v$; on a 
$$\left\{\begin{array}{ll}a_1=&\frac{N}{16(J_1-1728)^2}\\
a_2=&\frac{3}{2}\,{\frac {{v}^{2}J_{{1}}-J_{{1}}+1728}{J_{{1}}-1728}}\end{array}
\right.$$
$$\begin{array}{ll}{\rm avec} \;N=& 
14929920-17280\,J_{{1}}+5\,{J_{{1}}}^{2}+10368\,{v}^{2}
J_{{1}}-6\,{v}^{2}{J_{{1}}}^{2}\\&-13824\,{v}^{3}J_{{1}}+8\,{v}^{3}{J_
{{1}}}^{2}+9\,{v}^{4}{J_{{1}}}^{2}.
\end{array}$$

La courbe $E_1$ est tordue de $E(J_1)$ par $v$.

\medskip

L'invariant de $E_2$ est
$J_2=J_1\;U^3/V^3$, avec
$$\begin{array}{ll}U=&-  2985984+35831808\,v+{J_{{1}}}^{2}-3456\,J_{{1}}-93312\,{v}^
{2}J_{{1}}+54\,{v}^{2}{J_{{1}}}^{2}\\&+81\,{v}^{4}{J_{{1}}}^{2}
+108\,{
v}^{3}{J_{{1}}}^{2}+12\,v{J_{{1}}}^{2}-41472\,vJ_{{1}}-186624\,{v}
^{3}J_{{1}}-559872\,{v}^{4}J_{{1}} \\
&{\rm et}\\

V= &243\,{v}^{4}{J_{{1}}}^{2}+216\,{v}^{3}{J_{{1}}}^{2}+54\,{v}
^{2}{J_{{1}}}^{2}-{J_{{1}}}^{2}+3456\,J_{{1}}\\&-2985984-93312\,{v}^{2}J
_{{1}}-373248\,{v}^{3}J_{{1}} .\end{array}$$

La courbe $E_2$ a comme \'equation $y^2=x^3+Cx+D,$
avec $$\left\{\begin{array}{ll}
C=&-3v^2J_1(J_1-1728)J_2V^3\\D=&2v^3J_1(J_1-1728)(J_2-1728)V^3\end{array}
\right.$$

$6$) R\'eciproquement, \`a tout couple d'invariants distincts $J_1$ et $J_2$ dans
$k-\{0,1728\}$
correspondent douze courbes $C_{a_1,a_2}$ d\'efinies sur $\overline{k}$,
dont les jacobiennes sont isog\`enes \`a $E(J_1)^2\times E(J_2)$. Plus
pr\'ecis\'ement, soit $h$ l'une des $12$ homographies envoyant les abscisses des 
points d'ordre $3$ de
$E(J_1)$ sur les abscisses des points
d'ordre $3$ de $E(J_2)$;  
le param\`etre $v$ du point $5)$ est donn\'e par
$$\frac{1}{v}=3h(\infty).$$

}

\medskip
\ni
{\sc Remarque $2.1$.} 

\ni
Soit $\Phi_2=\Phi_1\circ \sigma:\;C_{a_1,a_2}\ra E_1$, et
$\Psi:\;C_{a_1,a_2}\ra E_2=C_{a_1,a_2}/<\sigma>$ un rev\^etement
de degr\'e $3$ tel que $\Psi(\rho:\rho^2:1)=0$.  
Soit $f:\;E_1^2\times E_2\ra \Pic^0(C_{a_1,a_2})$
d\'efinie par $$(P,Q,R)\ms \Phi_1^*((P)-(0))+\Phi_2^*((Q)-(0))+\Psi^*((R)-(0))$$
et $g:\;\Pic^0(C_{a_1,a_2})\ra E_1^2\times E_2$ d\'efinie
par $$D\ms (\Phi_1(D),\Phi_2(D),\Psi(D));$$    l'application compos\'ee
$g\circ f:\;E_1^2\times E_2\ra E_1^2\times E_2$ est donn\'ee par
la matrice $\left(\begin{array}{ccc}2&-1&0\\-1&2&0\\0&0&3\end{array}\right).$
Par suite, le degr\'e de $f$ (resp. de $g$)
est \'egal \`a $9$.

La jacobienne de $C_{a_1,a_2}$ est obtenue par recollement ([Serre Harvard], p. $37$ et seq. ou [Lauter2], p. $95$) 
de $E_1^2$
muni de la polarisation 
$\left(\begin{array}{cc}2&-1\\-1&2\end{array}\right)$ et de $E_2$ munie de la polarisation $3 \;Id$, via 
l'isomorphisme galoisien de $E_1[3]$ sur $E_2[3]$ du point $4)$ ci-dessus.

\subsection{Exemples}
On peut  faire une recherche syst\'ematique
des 
courbes optimales de genre $3$ de la forme $C_{a_1,a_2}$, par l'algorithme suivant: soit $k$ un corps
fini, de cardinal $q$, et 
$J$ l'ensemble, \'eventuellement vide, des invariants des 
courbes elliptiques dont la valeur absolue de la trace du Frobenius
sur $\F_q$ vaut $m_q$\footnote{On exclut les cas $m_q^2-4q=-3,-4,-8,-11$, car d'apr\`es [Lauter2], Appendice,
il n'existe pas de courbe optimale dans ces cas.}. Pour tout couple $j,j'$ d'\'el\'ements
de $J$, on calcule les racines de l'\'equation $J_2(j,v)=j'$,
puis le nombre de points de la courbe de genre $3$ associ\'ee \`a
chacune de ces racines par les formules ci-dessus. 

Pour les nombres premiers $p\leq 10723$, i.e. les $1308$ premiers
nombres premiers, cette m\'ethode \'echoue \`a trouver une courbe
de genre $3$ dont le nombre de points est $p+1+3m$ dans $97$ cas.

Pour ces cas, le nombre de classes de $m^2-4p$ est en
g\'en\'eral \'egal \`a $1$,  sauf dans $10$ cas, o\`u on trouve
$-20,-24,-35,-100,-123,-187$, de nombre de classes $2$, et $-59$,
de nombre de classes $3$.

Les $19$ nombres premiers $<1000$ pour lequel la m\'ethode pr\'ec\'edente
ne permet pas de trouver une courbe optimale sont
$53\;\footnote{Comme me l'a indiqu\'e C. Ritzenthaler, M. Nitzberg a montr\'e que $x^4+y^4+z^4=0$ est optimale pour $p=53$ 
([Serre Harvard], p. $72$).}
,167,173,193,293,311,347,$ 
$353,359,479,523,557,569,661,709,773,787,823,997.$

\medskip
On peut faire de m\^eme lorsque $q=p^n$ n'est pas premier, et
jouer alors avec le fait que $J_1$ et $J_2=J_1^p$ peuvent
\^etre distincts:
soit par exemple $q=19^3$; en prenant 
$J_1=J_2$, o\`u $J_1$ est l'invariant d'une courbe elliptique optimale
sur $\F_q$, on trouve une courbe de genre $3$ dont le nombre de points
est $q+1-3m_q$; par contre,
en prenant $J_2=J_1^{19}$, on trouve une courbe optimale.

\subsection{Cas o\`u $J_1=J_2$}
On cherche ici les courbes $C_{a_1,a_2}$ telles que $J_1=J_2$; 
avec les notations pr\'ec\'edentes, $J_1/J_2=M_3^3$,  donc $M_3$ est
alors \'egal \`a $1,\rho$ ou $\rho^2$.

\subsubsection{$M_3=1$}
L'\'equation $M_3=1$ s'\'ecrit
$$\begin{array}{l}-288-720\,a_{{2}}-936\,{a_{{2}}}^{2}+2592\,a_{{1}}a_{{2}}-648\,{a_{{2}
}}^{3}+1944\,a_{{1}}{a_{{2}}}^{2}\\-243\,{a_{{2}}}^{4}+1872\,a_{{1}}-
3888\,{a_{{1}}}^{2}=0\end{array}$$
courbe rationnelle que l'on peut param\'etrer par
$$\left\{
\begin{array}{l}
a_1={\frac {1}{432}}\,{\frac {112\,{t}^{4}+272\,{t}^{3}+408\,{t}^{2}+296\,
t+127}{ \left( {t}^{2}+t+1 \right) ^{2}}},\\
a_2=-{\frac {8\,{t}^{2}+10\,t+9}{6({t}^{2}+t+1)}}
\end{array} \right. $$

Dans ce qui suit, on note $C_t$ la courbe d'\'equation
$$(t^2+t+1)^2(a_1T_1^4+a_2T_1^2T_2+2T_1T_3+T_2^2)=0,$$ $a_1$ et $a_2$ \'etant les 
fractions rationnelles en $t$ ci-dessus.

On a alors $J_1=1728t^3$; la courbe $E_2$ a comme \'equation
$y^2=x^3-3t(t^3-1)x-2(t^3-1)^2$, 
et la courbe $E_1$ est tordue de $E_2$ par $-3(t^2+t+1)$.

Le discriminant de  $C_t$ est
$2^{118}3^{74} (t-1)^9(t^2+t+1)^{43}.$

\medskip
\ni
{\sc Remarques $2.3.1$.}

\medskip
1) D'apr\`es [Serre], p. $305$, l'invariant modulaire d'une courbe elliptique $E$
est un cube si et seulement si le degr\'e de l'extension obtenue par
adjonction des points d'ordre $3$ de la courbe est premier \`a $3$.

\medskip
2) Supposons que $k=\F_q$ ($q$ impair).

a)  Si $q\equiv 2 \mod 3$,  $j$ est toujours
un cube (comme tout \'el\'ement de $k$).

b) Si $q\equiv 1 \mod 3$, une condition suffisante pour que $j$ soit un
cube est que $a_q=q+1-N_q(E)$ soit divisible par $3$; en effet, dans ce
cas, ni $E$ ni sa tordue quadratique n'ont de point d'ordre $3$ rationnel
sur $\F_q$, et le polyn\^ome $f_3$ des abscisses des points d'ordre $3$
de $E$ et de sa tordue n'a donc pas de racine dans $\F_q$. Par suite,
il n'a pas de facteur irr\'eductible de degr\'e $3$.

Par ailleurs, 
la relation
$$(t^3-1)^2=(t^2+t+1)((t)^2+t+1)((t\rho^2)^2+t\rho^2+1)$$ 
montre que,
d\`es que $j$ est un cube, il existe $t \in k$ tel que
$j=1728 t^3$ et $t^2+t+1$ est un carr\'e dans $k$.

\medskip
3) Lorsque $t^2+t+1=0$, i.e. $t=\rho$ ou $\rho^2$, on a $J_1=J_2=1728$; la
quartique $C_{t}$ a un discriminant nul, et sa r\'eduction modulo $t^2+t+1$ est, sur
une extension convenable de $k(t)$, la courbe hyperelliptique $y^2=x(x^6-1)$. 
Pour le 
voir, il suffit d'\'ecrire l'\'equation de la quartique sous la forme $A^2+\eps B+O(\eps^2)=0$
([Clemens], p. $155-157$, ou encore [Elkies], p. $82$), $A=0$ \'etant l'\'equation d'une conique $Co$   non singuli\`ere;
lorsque $\eps$ tend vers $0$, la quartique "tend" vers la courbe hyperelliptique
rev\^etement double de la conique $Co$ ramifi\'e aux points d'intersection de $Co$ et
de la courbe d'\'equation $B=0$. 

Ici, on param\`etre la conique $t^2+t+1=u^2$ 
par $\left\{\begin{array}{l}
t=\rho^2\frac{\eps^2-\rho^2}{\eps^2-1}\\
u=-\sqrt{-3}\frac{\eps}{\eps^2-1}
\end{array}\right. .$

Pour $\eps=0$, on a $(t,u)=(\rho,0)$.

Apr\`es le changement de variables $(X,Y,Z)\ms
(X,Y,S)$, avec $X+Y+Z=2Su$,
la quartique $C_t$ s'\'ecrit 
$A^2-12\eps B+O(\eps^2)=0,$ 
o\`u $A=S^2-\sqrt{-3}(X^2+XY+Y^2)$ et $B=S(X+Y)(A-\sqrt{-3}XY).$

En param\'etrant 
la conique $Co$ d'\'equation $A=0$ par $\mu=\sqrt[4]{-3}(X-\rho Y)/S$, les
huit 
param\`etres des points
d'intersection de $Co$ et de la quartique d'\'equation $B=0$ sont   
$0,\infty,\pm \rho,
\pm \rho^2,
\pm 1$; par suite,
sur $\overline{k}$, la courbe hyperelliptique "limite" de $C_t$ quand
$t$ tend vers $\rho$ 
a comme \'equation $y^2=x(x^6-1)$.

\subsubsection{Cas o\`u $k=\Q$ et $E_1$ de type CM}
Si $k=\Q$, $E_1$ et $E_2$ ne sont pas $k$-isomorphes,
puisque tordues l'une de l'autre par $-3(t^2+t+1)$. N\'eanmoins, si elles
sont \`a multiplications complexes, elles peuvent \^etre isog\`enes sur $\Q$;
dans ce cas, la trace du Frobenius $\mod p$ de $C_t$, pour $p$ premier et premier
au discriminant de la courbe, est \'egale au triple de celle de $E_1$.

\medskip

Les $13$ courbes elliptiques \`a multiplications complexes d\'efinies sur 
$\Q$  ont comme invariant
$$2^6 3^3,2^6 5^3,0,-3^3 5^3,-2^{15},-2^{15} 3^3,-2^{18} 3^3 5^3,-2^{15} 3^3 5^3 11^3,$$ $$-2^{18} 
3^3 5^3 23^3 29^3,2^3 3^3 11^3,2^4 3^3 5^3,3^3 5^3 17^3,-3 \;2^{15} 5^3,$$
les anneaux d'endomorphismes associ\'es  \'etant ceux de discriminant
$$-4,-8,-3,-7,-11,-19,-43,-67,-163,-16,-12,-28,-27.$$

\`A part le onzi\`eme et le treizi\`eme, tous ces invariants sont des
cubes, les valeurs de $t$ associ\'ees \'etant
$$1,5/3,0,-5/4,-8/3,-8,-80,-440,-53360,11/2,{\frac {85}{4}}.$$

La valeur $t=1$ donne une quartique singuli\`ere; pour les dix autres
valeurs, on obtient dix courbes de genre $3$. \`A des facteurs carr\'es pr\`es,
les quantit\'es $-3(t^2+t+1)$ valent respectivement
$-3, -3, -7, -3, -19, -43, -67, -163, -1, -7$, les anneaux d'endomorphismes
ayant comme discriminant $$-8,-3,-7,-11,-19,-43,-67,-163,-16,-28.$$ 

\`A part le premier et quatri\`eme cas, les courbes $E_1$ et $E_2$ sont donc
$\Q$-isog\`enes, et on obtient ainsi huit quartiques d\'efinies sur
$\Q$ dont la jacobienne est $\Q$-isog\`ene au cube d'une courbe elliptique
\`a multiplications complexes, les anneaux d'endomorphismes
\'etant $-3,-7,-19,-43,-67,-163,-16,-28.$

\medskip

{\sc Exemples.-} 1) Pour $t=-5/4$, on trouve une quartique $\overline{\Q}$-isomorphe
\`a la courbe de Klein, d'\'equation $x^3y+y^3+x=0$, mais qui ne lui est pas $\Q$-isomorphe (la courbe
de Klein n'ayant pas d'involution d\'efinie sur $\Q$). Plus
pr\'ecis\'ement, l'isomorphisme entre la courbe de Klein 
et la quartique ci-dessus est donn\'e par
la matrice $$\left(\begin{array}{lll} u_1&u_2&u_3\\u_3&u_1&u_2\\u_2&u_3&u_1
\end{array}\right)$$ avec $u_1^4$ racine de $X^3+26X^2+13X+1,$ 
$u_2=\frac{1}{29}u_1(5u_1^8+126u_1^4-30)$ et 
$u_3=\frac{1}{29}u_1(17u_1^8+440u_1^4+159).$

Si $\delta$ est le d\'eterminant de cette matrice, on a 
$\delta^4=-14^4$.

La remarque $3$ du paragraphe pr\'ec\'edent permet de retrouver le fait
que, modulo $7$, la quartique de Klein se r\'eduit en 
une courbe hyperelliptique. (De m\^eme, les quartiques associ\'ees aux courbes \`a multiplications
complexes par l'anneau des entiers de $-19,-43,-67,-163$  se r\'eduisent 
modulo $19,43,67,163$ en une courbe hyperelliptique).

\medskip
2) Pour $t=-5/4$ (resp. $t=85/4$), qui correspondent aux discriminants $-7$ et $-28$,
 on obtient une quartique dont la jacobienne
est $\Q$-isog\`ene au cube de la courbe $y^2=x^3-35x+98$ (resp. 
$y^2=x^3-1933155x -1034488098$). Ces deux courbes sont $2$-isog\`enes
sur $\Q(\sqrt{57})$, mais pas sur $\Q$. Par suite, si $p$ est un nombre
premier tel que $\left(\frac{p}{57}\right)=-1$, les traces du Frobenius
$\mod p$ sont oppos\'ees d'une quartique \`a l'autre. En particulier,
si $p=4b^2+7$, avec $b\geq 4$ (ce qui assure que la partie enti\`ere
de $2\sqrt{p}$ est \'egale \`a $2b$), et  si $\left(\frac{p}{57}\right)=-1$,  
on obtient ainsi une quartique optimale et une quartique minimale sur
$\F_p$. Par exemple, les nombres premiers $<10000$   v\'erifiant
ces conditions sont
$$\{151, 263, 331, 491, 907, 1031, 1163, 1451, 1607, 2311, 4363, 5483, 5783\}.$$

\medskip

3) Dans les cas de discriminant $-19,-43,-67$ et $-163$, on 
obtient des mod\`eles sur $\Q$ de courbes trouv\'ees par Ritzenthaler dans 
[Ritzenthaler],
\`a savoir celles  dont la derni\`ere entr\'ee dans la Table $2$ de 
[{\it ibid.}] vaut
$2\times 6$.

\subsubsection{Le cas o\`u $k=\F_{7^n}$}
{\sc Proposition $2.3.3$.}
{\it Soit $n$ un entier impair, $q=7^n$, et $E$ une courbe elliptique d\'efinie sur $k=\F_q$, 
dont la trace du Frobenius $a$ est congrue \`a $9,15$ ou $18\mod 21$.
Il existe une courbe de genre $3$, d\'efinie sur $k$, dont
le nombre de points vaut $7^n+1-3a.$}

En effet, on a alors $a\equiv 0\mod 3$, donc d'apr\`es ce qui 
pr\'ec\`ede, l'invariant de $E$ est un cube $1728 t^3$, o\`u
l'on peut choisir $t$ tel que $t^2+t+1$ est un carr\'e dans $k$.

Il existe donc une courbe de genre $3$ d\'efinie sur $k$ dont la
jacobienne est isog\`ene \`a $E_2^3$, o\`u $E_2$ a 
comme \'equation
$y^2=x^3+Ax+B$, avec $A=-3t(t^3-1,B=-2(t^3-1)^2.$

Or, si $E$  a comme
\'equation $y^2=x^3+\alpha x+\beta$,  on sait ([Manin], th\'eor\`eme $1$)
que $a$ est congrue modulo $7$ \`a la norme
de son invariant de Hasse, soit $(3\beta)^{1+7+\ldots +7^{n-1}}.$
Par hypoth\`ese, c'est un carr\'e dans $\F_7$. 

Comme la trace du Frobenius de $E_2$ est aussi 
un carr\'e, puisque $3B=(t^3-1)^2$, $E_2$ est donc $k$-isomorphe \`a $E$,
d'o\`u le r\'esultat. 

\medskip

{\sc Corollaire $2.3.3$.} {\it Soit $k=\F_q$ un corps fini de caract\'eristique
$7$; 
il existe une courbe de genre $3$ d\'efinie sur
$\F_q$ dont le nombre de points est $\geq q+1+3(m_q-11)$.} 
\subsubsection{Le cas $M_3=\rho$}
On suppose ici que $\rho\in k$; 
la courbe d'\'equation $M_3=\rho$
admet une param\'etrisation rationnelle par un param\`etre $t$ tel que
$J_1=J_2={\frac { \left( -1+24\,t \right) ^{3}}{{t}^{3} \left( -1+27\,t
 \right) }}.$
On reconna\^\i t  l'invariant modulaire de la courbe $E_t$  d'\'equation
$y^2+xy+ty=x^3$, i.e. la courbe elliptique universelle ayant un point
d'ordre $3$ (cf. $2.1$). Par ailleurs, alors que, dans le cas $M_3=1$, $E_1$
et $E_2$  sont tordues quadratiques l'une de l'autre par $-3(t^2+t+1)$,
ici $E_1$ et $E_2$ sont $k(t)$-isomorphes, et tordues quadratiques
par $(\rho -1)(36t+\rho-1)$ de la courbe $E_t$.   

\medskip
\ni
Supposons que $k=\F_{7^n}$, $n$ impair, et prenons $\rho=2$.  Si $E$ est 
une courbe elliptique sur $k$ dont le nombre de points est $7^n+1-a$,
et que $a\equiv -1\mod 3$, le courbe $E$
a un point d'ordre $3$. Soit $F$ la courbe quotient de $E$ par le groupe engendr\'e par 
ce point. D'apr\`es $2.1$, il existe
$t\in k$ tel que $E=E_t$ et, $-3$ \'etant un carr\'e
dans $k$, $F=F_{-t+1/27}$. La courbe de genre $3$ associ\'ee \`a $t$ (resp. 
$-t+1/27$) est isog\`ene au cube de $E$ ou de sa tordue quadratique, selon que $t+1$
est un carr\'e ou non (resp. au cube de $F$ ou de sa tordue quadratique, selon que 
$36(1/27-t)+\rho-1=-t$ est un carr\'e ou non). Comme, \`a un carr\'e pr\`es, le discriminant
de $E_t$ est \'egal \`a $t(1+t)$, on en d\'eduit que, si $a$ est impair, i.e. si $E$ n'a pas 
de point d'ordre $2$ sur $k$, ou si $a\equiv 0\mod 4$,  dans lequel cas, quitte \`a remplacer $E$  par une courbe $2$-isog\`ene, $E$  a trois points d'ordre $2$, il existe un unique r\'esidu quadratique dans l'ensemble 
$\{-t,1+t\}$.
 
Donc, si $a\not \equiv -1\mod 3$ et $a\not\equiv 2\mod 4$, il existe une courbe de genre $3$ isog\`ene \`a $E^3$ et
une autre isog\`ene au cube de la tordue de $E$. En rempla\c cant $E$ par sa tordue
quadratique, on a le m\^eme r\'esultat si $a\equiv 1\mod 3$ et $a\not\equiv 2\mod 4$, d'o\`u la partie $2$
du th\'eor\`eme  
\'enonc\'e en introduction. Le point $2$ du corollaire s'ensuit (et peut \^etre pr\'ecis\'e:
\`a moins que $m_q\equiv 36,51,58$ ou $78\mod 84$, le d\'efaut est au plus $6$).

\section{Le cas de caract\'eristique $3$}
Dans ce paragraphe, les corps consid\'er\'es sont de caract\'eristique
$3$.
\subsection{Le cas $a_3 \neq 2$}
On consid\`ere comme pr\'ec\'edemment la famille de courbes
d'\'equation
$$a_1T_1^4+a_2T_1^2T_2+a_3T_1T_3+T_2^2.$$ 

Dans ce paragraphe, on suppose $a_3\neq 2$. Par une transformation lin\'eaire de la forme
$$X'=X+aT_1,Y'=Y+aT_1,Z'=Z+aT_1,$$ qui est inversible pour tout $a$,
on se ram\`ene \`a $a_2=0$ en prenant $a=-\frac{a_2}{1+a_3}.$

On note d\'esormais, dans ce qui suit, $C_{a_1,a_3}$ la courbe d'\'equation
$$a_1T_1^4+a_3T_1T_3+T_2^2.$$

Son discriminant est $2(a_3+1)^9a_1^3a_3^{12}$. 

La courbe $E_1$ quotient de $C$ par l'une quelconque des involutions naturelles
de $C$ a comme \'equation

$$v^2=u^3+a_3(a_3+1)u^2+2a_3^2a_1,$$
son invariant modulaire \'etant
$$J_1=\frac{a_3(a_3+1)^3}{a_1},$$
un morphisme du mod\`ele affine de $C$ obtenu en posant $x=X/Z,y=Y/Z$  vers $E_1$ \'etant donn\'e par
$$u=\frac{-a_3}{T_1},\;\;
v=\frac{T_2-a_3T_1}{T_1^2}.$$

La courbe $E_2$ quotient de $C$ par le groupe d'ordre $3$ engendr\'e 
par $(x,y,z)\ms (y,z,x)$ 
a comme \'equation
$$v^2=u^3+a_3(a_3+1)u^2+a_1a_3^5,$$
le morphisme de $C$ vers $E_2$ \'etant donn\'e par

$$\left\{\begin{array}{l}
u=2a_3^2 \frac{a_1UV+a_3}{(U+V)(2+a_1(U+V))}\\
v=(U-V)\frac{a_1a_3^4+u^2}{a_3^2}\end{array}\right.$$
avec $U=\frac{x^2y+y^2+x}{xy}$ et $V=\frac{y^2x+x^2+y}{xy}.$
L'invariant de $E_2$ est $$J_2=\frac{2(a_3+1)^3}{a_1a_3^2}$$
et donc $$\frac{J_1}{J_2}=-a_3^3.$$

Par ailleurs, les coefficients de $u^2$ dans les \'equations de $E_1$
et $E_2$ sont les m\^emes, et sont non nuls. Comme la trace du Frobenius
de $E_1$ (resp. $E_2$) est la norme de l'invariant de Hasse de $E_1$,
qui en caract\'eristique $3$ est le coefficient de $u^2$ dans
une \'equation de Weierstrass de la courbe, les deux courbes $E_1$
et $E_2$ ont des traces de Frobenius \'egales et non nulles
modulo $3$; elles ne peuvent donc pas \^etre oppos\'ees. 

Par suite, si $k=\F_{3^n}$, et si $J_1$ et $J_2$ sont deux \'el\'ements
distincts et non
nuls de $k$,
il existe une courbe de genre $3$ 
d\'efinie sur $k$, dans la famille
ci-dessus. dont la jacobienne est isog\`ene \`a $E_1^2\times E_2$, o\`u
$E_1$ (resp. $E_2$) a comme invariant $J_1$ (resp. $J_2$): le coefficient
$a_3$ est d\'etermin\'e par $J_1/J_2=-a_3^3$, d'o\`u 
$a_1=\frac{a_3(a_3+1)^3}{J_1}.$

Supposons d\'esormais $n$ impair $\geq 3$, et $q=3^n$. Pour tout $j\in k-\F_3$, on 
obtient ainsi une courbe $C$ associ\'ee au couple $(J_1=j,J_2=j^3)$,
le coefficient $a_3$ correspondant v\'erifiant $a_3^3=-J_1/J_2=-1/J_1^2$,
et une courbe $C'$ associ\'ee au couple $(J_2,J_1)$, le coefficient $a'_3$
correspondant  \'etant \'egal \`a $1/a_3$.

Le coefficient en $x^2$ de la courbe $E_1$ est $a_3(a_3+1)$, alors que
celui de la courbe $E'_1$ est $a'_3(1+a'_3)$. Le quotient de ces deux 
termes vaut $a_3^3=-1/J_1^2$, et n'est donc pas un carr\'e. Par suite,
les courbes $E_1$ et $E'_1$ sont tordues quadratiques l'une de l'autre,
et la jacobienne de $C'$ est la tordue quadratique de celle de $C$.

D'apr\`es [Lauter2], p. $92$ et p. $99$, il existe $a$ tel que
$a^2-4 \times 3^n=-8$ (resp. $-11$)  si et seulement si $n=1$ ou $n=3$ (resp.
$n=1$ ou $5$). 
Par suite, pour $n\geq 7$, pour tout entier $a$ non divisible par $3$
et tel que $|a|\leq 2\sqrt{q}$, o\`u $q=3^n$, il existe une courbe elliptique
d'invariant dans $\F_q-\F_3$ et dont la trace du Frobenius vaut $a$. 
 
On en d\'eduit la proposition suivante:

\medskip

\ni
{\sc Proposition $3.1$ .} {\it Soit $n$ un entier impair $\geq 7$,  et $q=3^n$.  
Si
$m_q$ n'est pas divisible par $3$, 
il existe deux courbes de genre $3$ d\'efinies sur
$\F_q$, l'une optimale,
i.e. dont le nombre de points vaut $q+1+3m_q$, l'autre 
dont le nombre de points vaut $q+1-3m_q$.}

\medskip

\ni
{\sc Remarques $3.1$.} 

1) Il existe une infinit\'e de $n$ impairs pour lesquels
$m_q\not\equiv 0\mod 3.$ 
En effet, $m_q\equiv 0\mod 3$ si et seulement si le $(n-1)/2$-i\`eme
terme de la mantisse du d\'eveloppement $3$-adique 
$1.1101120222201212202001\ldots$ de $2\sqrt{3}$ est nul.
Comme $2\sqrt{3}$ n'est pas un nombre rationnel, les termes de ce
d\'eveloppement ne peuvent pas \^etre tous nuls \`a partir d'un certain rang.

Dans le cas o\`u $m_q\equiv 0\mod 3$, (ce qui, exp\'erimentalement, pour $n$ impair, semble
arriver dans une proportion d'environ $1/3$), d'apr\`es ([Serre Harvard], p. $7$ ou [Lauter2], Table $1$, le d\'efaut
est au moins $2$, et s'il est $2$ les traces du Frobenius
sont $-(m_q+1-4\cos ^2 \frac{\pi}{7})$ et ses deux conjugu\'ees.

Ceci exige en particulier que la partie fractionnaire de $m_q$ soit
sup\'erieure \`a $1-4\cos^2 \frac{3\pi}{7}.$ Par exemple, pour $n$ impair
$\leq 400$, ceci arrive pour $n=15, 47, 53, 69, 159, 329, 349, 375, 383, 399$.

\medskip

2) Il n'existe qu'un nombre fini de $n$ impairs
pour lesquels $m_q$ soit la valeur absolue de la trace du Frobenius sur
$k$ d'une courbe d'invariant $\pm 1$. 
Plus g\'en\'eralement, soit $E$ une courbe elliptique d\'efinie sur un corps fini $\F_q$,
de Frobenius $\pi$; si $a_{q^n}=\pi^n+\overline{\pi}^n$, 
il existe une constante $c$ telle que, pour
$n$ impair $\geq 3$, on a $|2q^{n/2}-|a_{q^n}||\geq \frac{q^{n/2}}{n^c},$
quantit\'e qui tend vers l'infini avec $n$.
Ceci d\'ecoule 
par exemple du fait que, si $\beta$ est un nombre alg\'ebrique, ici $\pm
\frac{\pi^n}{q^{n/2}}$, il existe
une constante $c(\beta)$ telle que, pour tout entier $n\geq 2$ tel que $\beta^n
\neq 1$,  on a
$$|\beta^n-1|>n^{-c(\beta)}.$$ 
\medskip

En fait, dans le cas $j=1$, o\`u $\pi=\pm \frac{1+\sqrt{-11}}{2}$, 
(resp. $j=2$, pour lequel
$\pi=\pm 1+\sqrt{-2}$) 
il 
semble que le seul $n$ pour lequel il existe une courbe d'invariant $j$ 
optimale
sur $\F_{3^n}$ soit $n=5$ (resp. $n=3$). 

Dans ces deux cas, il n'existe pas de courbe de genre $3$ optimale
sur $\F_{3^n}$, puisqu'on est dans les cas o\`u il n'existe pas de
module ind\'ecomposable de rang $3$ (cf. [Lauter2], Appendice).

\medskip

3) Dans le cas o\`u $q=3^5$ (donc $m_q=31$), 
Lauter montre ([Lauter2], p. $14$, Prop. $6$) que
pour toute
courbe $C$ de genre $3$ sur $k=\F_q$, on a $|q+1-\Card C(k)|\leq 3m_q-3$,
et qu'il existe une courbe $C_0$ pour laquelle on a \'egalit\'e, le signe
de $q+1-\Card C_0(k)$ n'\'etant pas connu 
({\it ibid.}, 
Note $4.1.5$).

Les courbes de genre $3$ ci-dessus permettent en fait de montrer 
qu'il existe, et de construire effectivement,  une courbe $C$ de genre $3$ telle que $\Card C(k)=3^5+1+3m_q-3.$
En effet, l'invariant de la courbe elliptique optimale sur $k$, dont la
trace du Frobenius vaut $31$,  
est $J_1=1$; soit $z$ une racine de $z^5-z+1$; $z$ engendre $k^*$, 
et les invariants des courbes elliptiques dont la trace du Frobenius
vaut $31-3=28$ sont $z^{38}, z^{77}$ et leurs conjugu\'es. 

Prenons $J_2=z^{38\times 3}.$ D'apr\`es ce qui pr\'ec\`ede, il existe une courbe $C_{a_1,a_3}$
d\'efinie sur $k$, avec $a_3=-1/J_2$, dont la jacobienne est $k$-isog\`ene
\`a $E_1^2\times E_2$, $E_1$ et $E_2$ d'invariants respectifs $J_1$ et $J_2$; l'invariant de Hasse de chacune
de 
ces deux courbes est $a_3(a_3+1)$, 
c'est-\`a-dire, \`a un carr\'e pr\`es, $1-J_2$; sa norme vaut $2$,
et est donc congrue \`a $-28=-(m_q-3)\mod 3$. Par suite, la trace du Frobenius
de $E_1$  est $-31$ et celle de $E_2$ est $-28$; donc la courbe de genre
$3$ associ\'ee a un nombre de points \'egal \`a $3^5+1+3m_q-3.$

\subsection{Le cas $a_3=2$}
Dans le paragraphe pr\'ec\'edent, on a suppos\'e $a_3\neq 2$; si $a_3=-1$,
le discriminant de la courbe $a_1T_1^4+a_2T_1^2T_2-T_1T_3+T_2^2=0$ 
est $-a2^9$;
par chacune des trois involutions naturelles, elle a comme
quotient la courbe elliptique d'invariant nul et d'\'equation
$y^2=x^3+a_2x+2a_1+a_2^2,$ et la courbe elliptique quotient
par le sous-groupe d'ordre $3$ engendr\'e par $(X,Y,Z)\ms (Y,Z,X)$
est aussi d'invariant nul, et a comme \'equation
$$y^2=x^3+a_2a_1^2x+2a_1^4.$$

Supposons que $-a_2$ ne soit pas un carr\'e dans $\F_q$, o\`u $q=3^n$.
L'application $x\ms x^3+a_2x$ \'etant bijective, le nombre de points
de chacune des deux courbes pr\'ec\'edentes est \'egal \`a $q+1$.
Par suite, dans $\F_{3^{2n}}$, ces courbes sont optimales, d'o\`u des
courbes de genre $3$ sur $\F_{3^{2n}}$ non hyperelliptiques et optimales. 

\bigskip

\ni
{\sc R\'ef\'erences.}

\medskip
\ni
[Auer-Top] Roland Auer, Jaap Top, Some Genus 3 Curves with Many Points. ANTS 2002.

\medskip
\ni
[Ciani] E. Ciani, I Varii Tipi Possibili di Quartiche Piane pi ` 
u Volte Omologico-Armoniche,{\it  Rend. 
Circ. Mat. Palermo} 13 (1899), 347-373. 

\medskip
\ni
[Elkies], The Klein Quartic in Number Theory, The Eightfold Way, 
MSRI Publications, 
Volume 35, 1998.

\medskip
\ni
[Ibukiyama] Ibukiyama, T.: On rational points of curves of genus 3 over finite fields, {\it Tohoku Math J.}
45 (1993), 311-329.

\medskip
\ni
[La-Ritz] Lachaud Gilles, Ritzenthaler Christophe.
On a conjecture of Serre on abelian threefolds.
{\it Algebraic Geometry and its applications (Tahiti, 2007)}, 88Ð115. World Scientific, Singapore, 2008

\medskip
\ni
[Lauter1] K. Lauter,
Geometric methods for improving the upper bounds on the number of rational points on algebraic curves over finite fields. Lauter, Kristin, with an appendix in French by J.-P. Serre. {\it J. Algebraic Geom.} 10 (2001), no. 1, 19--36.

\ni
\medskip
[Lauter2] K. Lauter, 
 The maximum or minimum number of rational points on genus three curves over finite fields, with an Appendix by J-P. Serre, {\it Compositio Math.}  134 (2002) 87-111.

\medskip
\ni
[Manin] Ju. I. Manin. The Hasse-Witt matrix of an algebraic curve. 
{\it Trans. Amer. Math. Soc.}, 45:245-246, 1965 (English translation of a Russian original)
 
\medskip

\ni
[Na-Rit1] E. Nart, C. Ritzenthaler. 
Genus 3 curves with many involutions and application to maximal curves in characteristic 2, 	{\it arXiv}:0905.0546v1

\medskip
\ni
[Na-Rit2] E. Nart, C. Ritzenthaler, Jacobians in isogeny classes of supersingular abelian
threefolds in characteristic 2, {\it Finite Fields and their applications}, 14, (2008), 676-702.

\medskip

\ni
[Ritzenthaler], C. Rltzenthaler, Explicit Computations of Serre's obstruction
for genus $3$ curves and application to optimal curves, arXiv:0901.2920.

\medskip

\ni
[Serre Harvard] J.-P. Serre, Rational Points on curves over finite fields. Notes by F. Gouvea of lectures at
Harvard University, 1985.

\medskip
\ni
[Serre], J.-P. Serre, Propri\'et\'es galoisiennes des points d'ordre fini des courbes 
elliptiques.  
{\it Inventiones mathematicae}, volume 15; 259-331. 

\medskip
\ni
[Top] J. Top, Curves of genus 3 over small finite fields, {\it Indag. Math.} (N.S.) 14, (2003),
275-283.

\medskip

\ni
{\sc D\'epartement de Math\'ematiques, Universit\'e Paris $7$,  \\$175$ rue du Chevaleret, $75013$ Paris, France}

\ni
E-mail: mestre@math.jussieu.fr

\end{document}